\documentclass{article}
\usepackage{amsmath,amsfonts}
\usepackage[all]{xy}
\usepackage{amssymb,amsthm}
\usepackage{hyperref}
\usepackage{url}

\newcommand{\Set}{\mathbf{Set}}

\newcommand{\calC}{\ensuremath{\mathcal{C}}}
\newcommand{\calD}{\ensuremath{\mathcal{D}}}
\newcommand{\calE}{\ensuremath{\mathcal{E}}}

\newcommand{\calL}{\ensuremath{\mathcal{L}}}
\newcommand{\calN}{\ensuremath{\mathcal{N}}}
\newcommand{\calS}{\ensuremath{\mathcal{S}}}

\newcommand{\opCat}[1]{\ensuremath{#1^\mathrm{op}}}

\newcommand{\twopl}[2]{\langle #1, #2\rangle}

\newcommand{\Psh}[1]{\widehat{#1}}

\newcommand{\sk}{\mathrm{sk}}

\newcommand{\Sk}{\mathrm{Sk}}

\newcommand{\bfC}{\mathbf{C}}

\newcommand{\eqf}{\mathfrak{f}}

\newcommand{\rmB}{\mathrm{B}}

\DeclareMathOperator{\Eq}{\mathit{Eq}}

\theoremstyle{plain}
\newtheorem{theorem}{Theorem}[section]
\newtheorem{corollary}[theorem]{Corollary}
\newtheorem{proposition}[theorem]{Proposition}
\newtheorem{lemma}[theorem]{Lemma}

\theoremstyle{definition}
\newtheorem{definition}[theorem]{Definition}
\newtheorem{example}[theorem]{Example}
\newtheorem{remark}[theorem]{Remark}

\title{
Non-singular maps in toposes \\ with a local state classifier
}
\author{Mat\'\i as Menni\footnote{{\em Consejo Nacional de Investigaciones Cient\'\i ficas y T\'ecnicas} and {\em Centro de Matem\'atica de La Plata, Universidad Nacional de La Plata, Argentina}}}
\date{\today}

\begin{document}

\maketitle

\begin{abstract}
Recent progress on the question of the size of the class of connected and hyperconnected geometric morphisms from a given topos has led to the definition  of  {\em local state classifier}. We  discuss a historical precedent which leads to the notion  of {\em non-singular map} and we show that, for  a topos $\calE$ with a local state classifier, and each object $X$ therein, the domain of the  full subcategory  of ${\calE/X}$ consisting of non-singular maps over $X$ is a topos, and that the inclusion is the inverse image functor of a hyperconnected geometric morphism. 
The prospective geometric applications direct our attention to local state classifiers in toposes `of spaces'. We show that, at least in the pre-cohesive topos of reflexive graphs, the local state classifier, which is a colimit by definition, may be characterized as a limit; more specifically, as a variant of a subobject classifier.
%
\end{abstract}

\tableofcontents

\section{Introduction}

Part of the motivation of the present paper is historical so it is relevant to start with that, following a short recollection of some of the relevant terminology.
The reader is  assumed to be familiar with some topos theory including the notion of  geometric morphism.

 Let ${p : \calE \rightarrow \calS}$ be a geometric morphism. We say that $p$ is {\em surjective} if the inverse image functor ${p^* :\calS \rightarrow \calE}$ is faithful. We say that $p$ is {\em connected} if $p^*$ is fully faithful. Notice that if $p$ is connected then it is trivially surjective. Sometimes, for emphasis, $p$ is said to be `connected surjective', instead of just `connected'. Finally, $p$ is said to be  {\em hyperconnected} if it is connected and the counit of ${p^* \dashv p_*}$ is monic.
Details may be found in   \cite[A4.6]{elephant}.

One of the main 1982 results in  \cite{Rosenthal1982} is that hyperconnected geometric morphisms from a  Grothendieck topos $\calE$ are in correspondence with certain structures, called {\em quotient systems}, on any site of $\calE$. It follows that the number of such maps from $\calE$ is small.
Rosenthal states in the first paragraph of his paper that the work there ``is an outgrowth of a suggestion by Lawvere \cite{LawvereSUNY-Lectures-1975-77} about unifying certain constructions involving $G$-sets, where $G$ is a group or a monoid''.

The original  version of \cite{Lawvere2005}, where  reflexive graphs are used to illustrate the distinction between `categories of spaces' and `generalized spaces', appeared in 1986. Further illustrations appear in \cite{Lawvere1989} where, in passing, Lawvere comments that, in presheaf toposes, there is an objective way to measure singularities. For this, he constructs, in any presheaf topos, an object $\Eq$ and a map ${\sigma_X : X \rightarrow \Eq}$ to `measures singularities'  in each presheaf $X$. We will recall the details in Section~\ref{SecQuotation}.

In April 2009, the  Peripatetic Seminar on Sheaves and Logic  88 took place in Cambridge to honour  Hyland and Johnstone on the occasion of their 60th birthdays. There,  Lawvere gave a talk entitled ``Open problems in topos theory'', an account of which may be found in \cite{LawvereOpenProblems}. The first problem discussed there asks if there is a Grothendieck topos  $\calE$ for which the number of connected  geometric morphisms from $\calE$ is not small, or if, at the other extreme, such geometric morphisms could be parameterized internally, as subtoposes are.

More recently, a partial (positive) answer to Lawvere's first problem was found in \cite{Hora2024}. Hora introduces the notion of {\em local state classifier} and shows that all Grothendieck toposes have one. It easily follows that if $\calE$ has a local state classifier then the number of hyperconnected geometric morphisms from $\calE$ is small. Since connected  geometric morphisms from a Boolean topos are hyperconnected, it follows that, if $\calE$ is Boolean and has a local state classifier, the  number of connected  geometric morphisms from $\calE$ is small \cite[Corollary~{5.9}]{Hora2024}. (A full solution was announced more recently in  \cite{HoraKamio2024}, but it plays no role in the present paper, because local state classifiers play no role in the (negative) answer to the more general question.)

Hora observes that the weaker version of his result, obtained by strengthening the hypothesis to a Boolean Grothendieck topos, follows already from Rosenthal's result (plus the observation that connected  implies hyperconnected for maps with Boolean domain), but that the techniques used in the two proofs are quite different.

In order to briefly discuss the technique in \cite{Hora2024} we introduce some notation.
For a category $\calC$, we let ${\calC_m \rightarrow \calC}$ be the bijective-on-objects subcategory determined by the monomorphisms in $\calC$.

\begin{definition}\label{DefLocalStateClassifier} 
A {\em local state classifier} in a category $\calC$ is a colimit of the inclusion ${\calC_m \rightarrow \calC}$.
\end{definition}

We stress that the definition does not make any assumption on the size of the category $\calC$.

\begin{definition}\label{DefCoherentFamOfMonos}
A  family ${(\rmB, \beta) = (\beta_X : \rmB X \rightarrow X \mid X \in \calC)}$ of monomorphisms  (in a category $\calC$) is 
{\em coherent}  if, for every monic ${m : X \rightarrow Y}$ in $\calC$, the left-bottom composite below factors through the right vertical map 
\[\xymatrix{
\rmB X \ar[d]_-{\beta_X} \ar@{.>}[r] & \rmB Y \ar[d]^-{\beta_Y} \\
X \ar[r]_-m & Y 
}\]
and the resulting square is a pullback.
\end{definition}

In other words, a coherent family of monomorphims in the category $\calC$ is a functor ${\rmB : \calC_m \rightarrow \calC_m}$ together with a natural transformation $\beta$ from $\rmB$ to the identity on $\calC_m$ and such that every naturality square of $\beta$ is a pullback.

 For instance, it is well-known that the counit of any hyperconnected geometric morphism is a coherent family of monomorphisms.
  Incidentally, Hora works with coherent families of {\em subobjects}, but the difference is inessential. Coherent families of monomorphisms are more directly connected to natural transformations, and we can consider the associated family of subobjects when necessary. We will do so without warning.

\begin{proposition}[{\cite[Proposition~{4.5}]{Hora2024}}]\label{PropHora}
 If $\calE$ is a topos with a local state classifier $\Xi$  then the following three concepts correspond bijectively:
 \begin{enumerate}
 \item Subobjects of $\Xi$.
 \item Morphisms ${\Xi \rightarrow \Omega}$.
 \item Coherent families of subobjects.
 \end{enumerate}
\end{proposition}
\begin{proof}
The correspondence between the first two items is immediate. 
That between the last two follows because a family  ${(m_X : S_X \rightarrow X \mid X \in  \calE)}$ of subobjects is coherent if, and only if, the corresponding family ${(X \rightarrow \Omega \mid X \in \calE)}$ of characteristic maps is a cocone for ${\calE_m \rightarrow \calE}$.
\end{proof}

We reproduced \cite[Proposition~{4.5}]{Hora2024} so that the reader may immediately see  how  easily the colimit property of the local state classifier implies that the number of coherent families of subobjects coincides with the number of maps between two fixed objects.
The main technical result in \cite{Hora2024} is, I think, that such colimits abound; more specifically, as mentioned above, that every Grothendieck topos has a local state classifier. The result has two proofs, one in \cite[Proposition~{3.21}]{Hora2024} and the other, more general, in \cite[Proposition~{B2}]{Hora2024} which does not require an underlying topos.

The present paper was motivated by the desire to understand the existence of local state classifiers in toposes, the relation with Rosenthal's work, and the observation that the construction in \cite[Proposition~{3.21}]{Hora2024} coincides with that used by  Lawvere in 1989 in order to measure singularities.

Hora's result shows that local state classifiers are useful to give a  bound to the size of certain classes of geometric morphisms.
Lawvere's exposition suggests the possibility of geometric applications. In order to start the exploration of these, we isolate the definition of non-singular map in a topos with a local state classifier and we show that the category of non-singular maps with a fixed codomain $X$ is a topos.
We invite the reader to think of this topos as a `petit' topos of the object $X$.

In Section~\ref{SecShells} we briefly discuss a variant of the notion of coherent family of monos that we call {\em shells}.
Shells will help us illustrate coherent families and also emphasize certain aspects of our discourse.
Section~\ref{SecRosenthal} adapts   Rosenthal's work in order to describe shells on presheaf toposes in terms of certain structures in the corresponding sites.
Lawvere's 1989 paper is quoted in Section~\ref{SecQuotation} and his construction is related with Hora's.
Singular maps are defined in Section~\ref{SecNonSingularMaps} and the `petit' categories of singular maps (over a fixed object) are proved to be toposes in Section~\ref{SecPetitToposes}, where we also show that non-singular maps do not form a calibration in the sense of \cite{Johnstone2012} and where, for reasons that become evident there, we show that certain class of hyperconnected geometric morphisms is small, regardless of whether their domain has a local state classifier or not.  Following \cite{Lawvere2005} we use reflexive graphs to illustrate the  petit toposes of non-singular maps  in Section~\ref{SecExample}. As it was perhaps to be expected, local state classifiers have different properties depending on whether they lie in a topos `of spaces' or in a `generalized locale'. In Section~\ref{SecLSCinPrecohesiveToposes} we prove a suggestive characterization in the former case.

\section{Shells}
\label{SecShells}

In this section we discuss a class of comonads that will be useful to structure the discourse of the paper and also to illustrate the notion of coherent family of monomorphisms.

\begin{definition}\label{DefShell}
A {\em shell} is an idempotent comonad $(\rmB,\beta)$ whose   counit $\beta$  is (pointwise) monic, and also mono-cartesian in the sense that, for every monic ${m : X \rightarrow Y}$, the naturality square
\[\xymatrix{
\rmB X \ar[d]_-{\beta_X} \ar[r]^-{\rmB m} & \rmB Y \ar[d]^-{\beta_Y} \\
X \ar[r]_-m & Y 
}\]
is a pullback.
\end{definition}

In other words, a shell is an idempotent comonad whose counit is a coherent family of monomorphisms in the sense of Definition~\ref{DefCoherentFamOfMonos}.

For instance, the comonad determined by a  hyperconnected geometric morphism is a shell.
 This is implicit in \cite[A4.6.6]{elephant}, more explicitly in \cite[Lemma~{2.6}]{Menni2021}, and also in \cite[Lemma~{4.4}]{Hora2024}.
If one thinks of such a geometric morphism ${p : \calE \rightarrow \calS}$ as a topos `of spaces' lying over a topos `of sets' $\calS$ then, for any space $X$, the counit ${p^* (p_* X) \rightarrow X}$ is pictured as the discrete subobject of points of $X$ or, roughly speaking, as the discrete skeleton of $X$.

\begin{definition}\label{DefCartesian} A map ${f  : X \rightarrow Y}$ is {\em cartesian} with respect to a comonad ${(\rmB, \beta)}$ if the  square
\[\xymatrix{
\rmB X \ar[d]_-{\beta_X} \ar[r]^-{\rmB f} & \rmB Y \ar[d]^-{\beta_Y} \\
X \ar[r]_-f & Y 
}\]
is a pullback.
\end{definition}

Of course,  monomorphisms are cartesian with respect to any shell.

If ${p : \calE \rightarrow \calS}$ is a hyperconnected geometric morphism, and ${f : X \rightarrow Y}$ is a map  in $\calE$, the subobject of $X$ in the pullback on the left below
\[\xymatrix{
\mathfrak{F} f \ar[d] \ar[r] & p^*(p_* Y) \ar[d]^-{\beta_Y} && p^*(p_* X) \ar[d]_-{\beta_X} \ar[rr]^-{p^* (p_* f)} &&   p^*(p_* Y) \ar[d]^-{\beta_Y} \\
X \ar[r]_-f                            & Y                                                && X \ar[rr]_-f                                                                    && Y
}\]
is called the subobject {\em of fibers} of $f$, and $f$ is said to have {\em discrete fibers} if the square on the right above is a pullback, that is, if $f$ is cartesian with respect to $\beta$. 
Recall that discrete fibers is a typical feature of local homeomorphisms (in the classical sense among topological spaces). Partly for this reason,  \cite{Lawvere2005} shows that, for any object $Y$ in $\calE$, the full subcategory ${\calS(Y) \rightarrow \calE/Y}$ of maps  (over $Y$)  with discrete fibers is the inverse image of a hyperconnected geometric morphism, presenting the topos ${\calS(Y)}$ as a candidate petit topos of `pseudo-classical sheaves' over  the space $Y$. These ideas are illustrated using the pre-cohesive ${p : \Psh{\Delta_1} \rightarrow \Set}$ of reflexive graphs over sets. In particular, \cite[Proposition~3]{Lawvere2005} states that, for the reflexive graph $L$ consisting of a single node and a single non-identity loop, the topos ${\calS(L)}$ coincides with the étendue of non-reflexive graphs.  See  also \cite{Menni2022}.

The shells determined by hyperconnected maps have a finite-limite preserving underlying functor.
We next discuss a different sort of example.

Let ${j : \calE_j \rightarrow \calE}$ be a level of $\calE$, that is,  an essential subtopos. The left adjoint to $p^*$ will be denoted by ${j_!}$ as  usual.
Lawvere's Dimension Theory  invites us to think of the leftmost adjoint ${j_! : \calE_j \rightarrow \calE}$ as the full subcategory of objects $X$ such that ${\dim X \leq j}$. Although it is not so in general,  the comonad ${j_!  j^*}$ is a shell in many cases of interest. For instance, every level of the topos of simplicial sets satisfies this property. The counit ${j_! (j^* X) \rightarrow X}$ is thought of as the $j$-(dimensional-)skeleton of $X$. With this intuition in mind, the pullbacks in Definition~\ref{DefShell} say that the $j$-skeleton  of a subobject is just the restriction of the $j$-skeleton  of the codomain.

More generally, if we take the images of the components of the counit of ${j_! \dashv j^*}$, we obtain a family of monomorphisms that is also natural.
Sometimes, the resulting {\em principal} comonad is a shell. For example, if $\calD$ is a small category where every map factors as a split epi followed by a monic, then the levels of the topos  ${\Psh{\calD}}$ whose associated principal comonad is a shell are in bijective correspondence  with the full subcategories of $\calD$ that are closed under subobjects \cite[Proposition~{3.2}]{Menni2024}.

So, a shell is something like a skeleton, but not quite. 

A related intuition is that, for a shell ${(\rmB, \beta)}$, ${\beta_X : \rmB X \rightarrow X}$ is the subobject of `singular figures of type ${(\rmB, \beta)}$'.
`Singular', in this vague sense, includes the idea of a figure that is below a certain dimension, as the examples with levels suggest.

Since shells are particular cases of coherent families of monomorphisms, \cite{Hora2024}  implies that the number of shells on a Grothendieck topos is small and, 
adapting some of the results in \cite{Rosenthal1982}, we will give a concrete description of shells on presheaf toposes.

\section{Shells on presheaf toposes}
\label{SecRosenthal}

\newcommand{\prP}{\mathrm{P}} 
\newcommand{\eqe}{\mathfrak{e}} 

The notion of an equivalence relation in a regular category is well-known; see, for example, \cite[A1.3.6]{elephant} and the references therein.
In particular, one may consider equivalence relations on representable objects in presheaf toposes, and these may be described in terms of the site as done in \cite{Lawvere1989}.

Let $\bfC$ be a small category. An {\em  equivalence relation} on an object $C$ in $\bfC$ is
meant the specification for every ${D\in \bfC}$ of a set of ordered pairs ${D \rightrightarrows C}$
of morphisms in $\bfC$ which is reflexive, symmetric, and transitive for each $D$
and which is closed with respect to composition by arbitrary  ${D'\rightarrow D}$.

Notice that an equivalence relation on $\bfC$ as defined in the previous paragraph is essentially the same thing as an equivalence relation, in the sense of the first paragraph,  on the representable ${\bfC(-,C)}$ in the topos ${\Psh{\bfC}}$ of presheaves on $\bfC$. 
We will confuse the two things.
For instance, an equivalance relation on $C$ in either sense, has a coequalizer ${\bfC(-, C) \rightarrow Q}$ in ${\Psh{\bfC}}$ which, by exactness, has the given equivalence relation as kernel pair.
Similarly, any map ${\bfC(-,C)\rightarrow X}$ in $\Psh{\bfC}$ has a kernel pair which we may identify with an equivalence relation on $C$.

If $\eqe$ is an equivalence relation on $C$ and ${f  : B \rightarrow C}$ is a map in $\bfC$ then the pairs ${g_1, g_2}$ of parallel maps with codomain $B$ such that ${f g_1 , f g_2 }$ are in $\eqe$ form an equivalence relation on $B$ that we denote by ${ \eqe \cdot f}$.

\begin{definition}\label{DefProbe} 
A {\em probe} in $\bfC$ is function $\prP$ that assigns to each $C$ in $\bfC$ a set ${\prP C}$ of equivalence relations on $C$ such that, for every ${f : B \rightarrow C}$ in $\bfC$ and $\eqe$ in ${\prP C}$, ${\eqe\cdot f}$ is in ${\prP B}$. 
\end{definition}

The class of probes introduced below will play a relevant role.

\begin{definition}\label{DefSaturated}
We say that a probe $\prP$ in $\bfC$ is {\em saturated} if, for every equivalence relation $\eqe$ on $C$ containing one in $\prP C$, ${\eqe \in \prP C}$.
\end{definition}

It should be stressed that the notion  of probe is a weak version of Rosenthal's  {\em system of quotients} which, as \cite[Remark (1), p.430]{Rosenthal1982} notes, are of a nature dual to Grothendieck topologies. While Grothendieck topologies select subobjects of representables, probes select quotients of representables. It is natural to wonder if this intuitive observation may be strengthened to an actual duality. 
Fragments of such a thing may be perceived in the relation discussed in \cite{Menni2021} between principal topologies in a topos  and certain hyperconnected quotients of the same topos; see also Proposition~\ref{PropGoodHyperconnectedAreFew}.

Let ${(\rmB, \beta)}$ be a coherent family of monomorphisms in $\Psh{\bfC}$.
For each $C$ in $\bfC$, let ${\prP_\beta C}$ be the set of equivalence relations on $C$ whose corresponding quotients ${\bfC(-,C) \rightarrow Q}$ are such that   $\beta_Q$ is an isomorphism.

\begin{lemma}\label{LemShellsInduceProbesNew} The function $\prP_\beta$   is a  probe in $\bfC$. If ${(\rmB, \beta)}$ is shell then $\prP_\beta$ is saturated.
\end{lemma}
\begin{proof}
To prove that $\prP_\beta$ is a probe, let ${f : B \rightarrow C}$ be a map  in $\bfC$ and let ${\eqe\in P C}$.
Then the  map
\[ R = \bfC(-,B)/(\eqe\cdot f) \rightarrow \bfC(-,C)/\eqe = Q \]
 induced by $f$ is monic.
Since, ${\beta_Q}$ is an isomorphism by hypothesis, then so is ${\beta_R}$, because of the `coherence' condition.

Assume now that ${(\rmB, \beta)}$ is shell. To show that $\prP_\beta$ is saturated just recall that, for idempotent comonads with monic counit, coalgebras are closed under quotients. See, e.g., \cite[Lemma~{2.2}]{Menni2021}.
\end{proof}

The intuition is that each probe in $\bfC$ determines a choice of  generic `singular' objects  in $\Psh{\bfC}$ that may be used to indentify singular figures in arbitrary presheaves. More precisely:

\begin{definition}\label{DefPfigureNew} 
For a probe $\prP$ in $\bfC$,  a figure  ${x :  \bfC(-, C) \rightarrow X}$ in $\Psh{\bfC}$   is {\em ($\prP$-)singular} if its  kernel is  in  ${\prP C}$.
\end{definition}

Singular figures and composition may be immediately related as follows.

\begin{lemma}\label{LemMapsPreserveSingularitiesNew} 
If $\prP$  is a probe in $\bfC$ and $X$ is a presheaf on $\bfC$ then,  for every $C$ in $\bfC$ and ${x \in X C}$,  the following hold:
\begin{enumerate}
\item If ${f : B \rightarrow C}$ is a map in $\bfC$ and ${x \in X C}$ is $\prP$-singular then so is ${x\cdot f}$.
\item For every ${\psi : X \rightarrow Y}$ in $\Psh{\bfC}$: 
\begin{enumerate}
\item If $\prP$ is saturated and  ${x \in X C}$ is $\prP$-singular then ${\psi_C x \in Y C}$ is $\prP$-singular.
\item If ${\psi}$ is monic  then,   ${x \in X C}$  is $\prP$-singular if, and only if, ${\psi_C x \in Y C}$ is.
\end{enumerate}
\end{enumerate}
\end{lemma}
\begin{proof}
In the commutative diagram below
\[\xymatrix{
 \bfC(-,B)  \ar@{->>}[d] \ar[rr]^-{\bfC(-,f)} && \bfC(-,C) \ar@{->>}[d]  \ar[r]^-x & X \\
 R \ar@{^{(}->}[rr] &&  Q \ar@{^{(}->}[ru]
}\]
the decorated arrows denote the evident epi/mono factorizations. The vertical figure  ${\bfC(-,C) \rightarrow Q}$ is the quotient of the kernel $\eqe$ of $x$, which is in ${\prP C}$ by hypothesis. By regularity, ${\bfC(-, B) \rightarrow R}$ is the quotient of ${\eqe\cdot f}$, so it is in ${\prP B}$ because $\prP$ is probe.
Also, the left-bottom composite shows that ${\bfC(-, B) \rightarrow R}$ is the quotient of the kernel of top composite ${\bfC(-,B) \rightarrow X}$.
By exactness the kernel must coincide with ${\eqe\cdot f  \in \prP B}$. So the composite figure is $\prP$-singular.

The kernel of ${\psi x}$  contains the kernel of $x$.
So, if $x$ is singular and $\prP$ is saturated then ${\psi x}$ is singular.

If $\psi$ is monic, the kernel of the composite
\[\xymatrix{
\bfC(-,C) \ar[r]^-x & \ar[r] X \ar[r]^-{\psi} & Y
}\]
coincides with the kernel of $x$.
\end{proof}

The concept introduced below is the analogue of   {\em $\calD$-generated object} (for a system of quotients $\calD$) formulated in  \cite[Definition~{1.4}]{Rosenthal1982}.

\begin{definition}\label{DefPskeletal}
For a probe $\prP$ in $\bfC$, we say that a presheaf  $X$ on $\bfC$  is {\em $\prP$-skeletal} if every figure of $X$ is  $\prP$-singular.
\end{definition}

For an arbitrary $X$ in $\Psh{\bfC}$ and $C$ in $\bfC$ we let 
\[\sk_{\prP, X, C} : (\Sk_\prP X) C   \rightarrow  X C \]
be the subset of $P$-singular figures of sort $C$.

\begin{lemma}\label{LemP-FCM} 
The assignment ${C \mapsto (\Sk_\prP X)C}$ extends to a presheaf on $\bfC$ that we denote  by  ${\Sk_\prP X}$. 
Also,  ${C \mapsto \sk_{\prP, X, C}}$ extends to a monic natural transformation ${\sk_{\prP, X} : \Sk_\prP X \rightarrow X}$ which is an isomorphism if, and only if, $X$ is $\prP$-skeletal.
Moreover, the family ${(\sk_{\prP, X} \mid X\in \Psh{\bfC})}$ is coherent.
\end{lemma}
\begin{proof}
The fact that  ${\Sk_\prP X}$ is a presheaf and  ${\sk_{\prP, X} : \Sk_\prP X \rightarrow X}$ is a natural transformation follows from the first item of Lemma~\ref{LemMapsPreserveSingularitiesNew}. (Notice that this is, essentially, the argument in \cite[p.~{427}]{Rosenthal1982}.)
It then easily follows that ${\sk_{\prP, X}}$  is an isomorphism if, and only if, $X$ is $\prP$-skeletal.

The family  ${(\sk_{\prP, X} \mid X\in \Psh{\bfC})}$ is coherent by the relevant item of Lemma~\ref{LemMapsPreserveSingularitiesNew}.
\end{proof}

To above lemma has an expected companion.

\begin{lemma}\label{LemPskeletonNew} If $\prP$ is saturated then ${(\Sk_\prP, \sk_{\prP})}$ is a shell, and its coalgebras are the $\prP$-skeletal objects.
\end{lemma}
\begin{proof}
If $\prP$ is saturated then  Lemma~\ref{LemMapsPreserveSingularitiesNew} implies that, for every ${\psi : X \rightarrow Y}$ in $\Psh{\bfC}$, the left-bottom composite below
\[\xymatrix{
\Sk_\prP X  \ar[d]_-{\sk_{\prP, X}} \ar@{.>}[r]^-{\Sk_\prP \psi} & \Sk_\prP Y \ar[d]^-{\sk_{\prP, Y}} \\
X \ar[r]_-{\psi} & Y
}\]
factors through the right vertical map above via a necessarily unique map that we may denote by ${\Sk_\prP \psi}$.
It easily follows that the assignment ${\psi \mapsto \Sk_\prP \psi}$ is functorial, so ${(\Sk_\prP, \sk_{\prP})}$ is a shell by Lemma~\ref{LemP-FCM}, which also implies that the coalgebras coincide with the $\prP$-skeletal objects.
\end{proof}

The rest of the section is devoted to the proof that the assignments 
\[ \prP \mapsto (\Sk_\prP, \sk_\prP)  \quad  \textnormal{ and } \quad (\rmB, \beta) \mapsto \prP_\beta \]
 form a bijective correspondence between shells on $\Psh{\bfC}$ and saturated probes in $\bfC$.

\begin{lemma}\label{LemAshellIsTheShellOfItsProbe} 
If  ${(\rmB, \beta)}$ is a coherent family of monomorphisms in  $\Psh{\bfC}$ then  ${(\rmB, \beta)}$ coincides with ${(\Sk_{(\prP_\beta)}, \sk_{(\prP_\beta)})}$. 
\end{lemma}
\begin{proof}
To lighten up the notation we let ${\prP = \prP_\beta}$.
To prove that ${\beta = \sk_\prP}$ let $C$ in $\bfC$.
Then, ${x \in (\Sk_\prP X) C \subseteq X C}$ if, and only if, $x$ is $\prP$-singular.
That is, if the kernel of ${x : \bfC(-, C) \rightarrow X}$ is in ${\prP C}$. 
Equivalently, if for the epi-monic factorization 
\[\xymatrix{
\bfC(-,C) \ar[d]_-q  \ar[rd]^-x  \\
Q \ar[r]_-{x'} & X
}\]
of $x$, ${\beta_Q}$ is an isomorphism.
This holds if, and only if,  $x'$ factors through $\beta_X$, because $x'$ is monic and $\beta$ is  coherent. 
Equivalently, $x$ factors through ${\beta_X}$.
\end{proof}

\begin{lemma}\label{LemAprobeIsAprobeOfItsShell} 
For any probe $\prP$,   ${\prP_{(\sk_{\prP})} = \prP}$.
\end{lemma}
\begin{proof}
An equivalence relation ${\eqe}$ on $C$ is in ${\prP_{(\sk_\prP)} C}$ if, and only if, the corresponding quotient ${q : \bfC(-,C) \rightarrow Q}$ is such that ${\sk_{\prP, Q}}$ is an isomorphism. Equivalently,  $q$ is a $\prP$-singular figure; which means that its kernel is in ${\prP C}$.
\end{proof}

To summarize, we have the following analogue of the results in \cite[Section~2]{Rosenthal1982}, which exhibit a correspondence between hyperconnected geometric morphisms from a Grothendieck topos $\calE$, and quotient systems on any site for $\calE$.

\begin{theorem}\label{PropShellsAndProbes} 
For any small category $\bfC$, the assignments ${\prP \mapsto (\Sk_\prP, \sk_\prP)}$ and ${(\rmB, \beta) \mapsto \prP_{\beta}}$ determine a bijection between  probes in $\bfC$ and coherent families of monomorphisms in ${\Psh{\bfC}}$, which ristricts to a bijection between saturated probes and shells.
\end{theorem}
\begin{proof}
The first part of the statement follows from Lemmas~\ref{LemAshellIsTheShellOfItsProbe} and \ref{LemAprobeIsAprobeOfItsShell}.
The restriction follows because, if $\prP$ is saturated then $\Sk_\prP$ is a shell by Lemma~\ref{LemPskeletonNew} and, if ${(\rmB,\beta)}$ is a shell then ${\prP_\beta}$ is saturated by Lemma~\ref{LemShellsInduceProbesNew}.
 \end{proof}

As suggested in \cite[p.~430, Remark~(3)]{Rosenthal1982}, if the category $\bfC$ came equipped with a Grothendieck topology $J$ then we could require a probe in $\bfC$ to interact well with $J$.  For instance, we could require that the equivalence relations are $J$-closed, but we will not pursue the idea here.

We mention in passing that, if ${(\prP_i \mid i \in I)}$ is  a family of probes in $\bfC$ then, if we let ${\prP C := \cap_{i\in I} \prP_i C}$,  it is easy to prove that the assignment ${C \mapsto \prP C}$  defines a probe $\prP$ in $\bfC$. Roughly speaking, probes are closed under intersection.
So given any family of equivalence relations in $\bfC$, we may consider the least probe that contains them.

Finally, it is relevant to mention that any probe in $\bfC$ determines a shell.
Indeed, if $\prP$ is a probe in $\bfC$ then, for any $C$ in $\bfC$, let ${\overline{\prP} C}$ be the set of equivalence relations on $C$ that contain an equivalence relation in ${\prP C}$.

\begin{lemma}\label{LemSaturation} 
The function ${C \mapsto \overline{\prP}C}$ is a saturated probe. It is the least saturated probe containing $\prP$.
\end{lemma}
\begin{proof}
To prove that $\overline{\prP}$ is a probe just notice that if ${\eqe \subseteq \eqf}$ are equivalence relations on $C$ with ${\eqe \in \prP C}$,  and ${f : B \rightarrow C}$ is a map in $\bfC$ then ${f^* \eqe \in \prP B}$ because $\prP$ is a probe and ${f^* \eqe \subseteq f^* \eqf}$.
\end{proof}

Of course,  a probe $\prP$ is   saturated if, and only if, ${\overline{\prP} = \prP}$.

\section{Measuring singularities}
\label{SecQuotation}

Let $\bfC$ be a small category.

In the pages before \cite[p.~282]{Lawvere1989} there are two ocurrences of the word ``singular''.
The first one is in p.~{265}:
\begin{quotation}
In case the objects and morphisms of $\bfC$ have some kind of geometrical interpretation,
it is often helpful to imagine that the more general objects of $\calS^{\opCat{\bfC}}$
push that interpretation to a natural limit: an object $A$ of $\bfC$ may be considered
in $\calS^{\opCat{\bfC}}$ as a generic ``figure'' and any ${\xymatrix{A \ar[r]^-x & X}}$ as a particular figure in $X$
(quite possibly singular, i.e. not necessarily monomorphic) of sort $A$.
\end{quotation}
and the second one in p.~{269}:
\begin{quotation}
\noindent ... it may be useful to consider quotients (like $C$) of objects of $\bfC$ (even
when ${C\not\in \bfC}$) as further ``generic'' figures, so that various types of singular
figures (such as loops) also become ``representable'' at least by objects of $\calS^{\opCat{\bfC}}$
\end{quotation}
which is illustrated with examples in reflexive graphs. Later, in p.~{282},  we read:
\begin{quotation}
As I have mentioned ``singularity'' several times, it may have occured to
the reader that there is an objective way of measuring it. Indeed for any
small category $\mathbf{C}$, there is a distinguished object  $\Eq = \Eq_{\mathbf{C}}$  in the topos
${\calS^{\opCat{\mathbf{C}}}}$  such that for every object $X$ there is a canonical map ${\xymatrix{X \ar[r]^-{\sigma_X} & \Eq}}$ in ${\calS^{\opCat{\mathbf{C}}}}$
 which does this. Namely the elements of sort $C$ of $\Eq$ are just all
the equivalence relations on $C$, where by an equivalence relation on $C$ is
meant the specification for every ${D\in \mathbf{C}}$ of a set of ordered pairs ${D \rightrightarrows C}$
of morphisms in $\mathbf{C}$ which is reflexive, symmetric, and transitive for each $D$
and which is closed with respect to composition by arbitrary  ${D'\rightarrow D}$. If
${C'\rightarrow C}$ in $\mathbf{C}$ and ${E\in\Eq(C)}$, then ${E\cdot\lambda \in \Eq(C')}$  is defined by taking for
each $D$
\[ \twopl{t_1}{t_2} \in E\cdot \lambda \Leftrightarrow \twopl{\lambda t_1}{\lambda t_2} \in E \]
thus making $\Eq$  into an object of ${\calS^{\opCat{\mathbf{C}}}}$. Of course it is more than just an
object, having a natural intersection operation ${\Eq \times \Eq \rightarrow \Eq}$ and greatest
point ${1\rightarrow \Eq}$  making it into a semilattice object, so in particular into an
ordered object. On the other hand, although the equality relation ${\Delta_C \in \Eq(C)}$ 
for each $C$ this is not natural, i.e. does not define a point ${\xymatrix{1 \ar[r]^-{\Delta} & \Eq}}$ {\em unless}
it happens that all morphisms in $\mathbf{C}$ are monomorphisms. The singularity
measurement ${\xymatrix{X\ar[r]^-{\sigma_X} & \Eq}}$ is defined, for each ${\xymatrix{C \ar[r]^-{x} & X}}$, by
\[ \sigma(x) = \{\twopl{t_1}{t_2} \mid x t_1 = x t_2 \} \]
the ``self-incidence''; then for ${\xymatrix{C' \ar[r]^-{\lambda} & C}}$ we  have
\[ \sigma ( x\lambda) = \sigma(x) \lambda \]
since ${(x\lambda) t_1 = (x\lambda) t_2}$  iff ${x(\lambda t_1)  = x(\lambda t_2)}$. The maps $\sigma_X$, although canonical,
are not natural when $X$ is varied, that is
\[\xymatrix{
X \ar[rd]_-{\sigma_X} \ar[rr]^-f & & \ar[ld]^-{\sigma_Y} Y \\
 & \Eq 
}\]
only commutes for ``non-singular'' $f$; of course
\[ \sigma_X \subseteq \sigma_Y \circ f \]
for all $f$, so we could say that $\sigma$ is natural in a suitable ``2-categorical'' sense.
Some $f$ might be ``equisingular''  in the sense that there exists an endomorphism $\lvert f \rvert$ of $\Eq$ so that a square commutes. 
For example with ${\mathbf{C} = \Delta_1}$, ${\Eq_{\mathbf{C}}}$ is the loop, and indeed in our generalization ${\mathbf{M}(T)}$ of ${\Delta_1}$, $\Eq$ itself is very singular.
\end{quotation}

Notice  that a probe in $\bfC$ is essentially the same thing as a subobject of $\Eq$, and that the probe is saturated if, and only if, the subobject is upper closed.

It is not unusual that the same original  idea sparks in different minds.
It is often considered to be a sign that the idea is good, especially if we restrict to scientific ideas.

\begin{theorem}[Hora]\label{ThmHoraExistenceOfLSC} The family ${(\sigma_X : X \rightarrow \Eq \mid X \in \Psh{\bfC})}$ is a local state classifier.
\end{theorem}
\begin{proof}
Just compare the definition of $\Eq$ in the quotation above with that of $\Xi$ in \cite[Example~{3.22}]{Hora2024}.
\end{proof}

\section{Non-singular maps}
\label{SecNonSingularMaps}

\newcommand{\aEq}{\mathfrak{S}}  

For general reasons, any meet semilattice ${(M,  \wedge, \top)}$ in a category $\calE$ with finite limits determines, for each object $X$, a partial order on the set of maps ${X \rightarrow M}$. In this restricted context, a family ${\mu = (\mu_X : X \rightarrow M \mid X \in \calE)}$ of maps will be called a    {\em lax cocone} if, for every ${f : X \rightarrow Y}$ in $\calE$,  ${\mu_X  \leq \mu_Y f}$ in the poset of maps from $X$ to $M$.

\begin{definition}\label{DefNonSingular}
A map ${f : X \rightarrow Y}$ is {\em non-singular (w.r.t. a family  ${\mu}$)} if ${\mu_X = \mu_Y f}$.
\end{definition}

\begin{remark}\label{RemNonSingularMonosAndLDcocones} A family  ${(\mu_X : X \rightarrow M \mid X \in \calE)}$ is a cocone for ${\calE_m \rightarrow \calE}$ if, and only if, every monomorphism is non-singular w.r.t. $\mu$.
\end{remark}

\begin{example}\label{ExNonSingularNaturalTransformations} With the notation of Section~\ref{SecQuotation},  a morphism ${f : X \rightarrow Y}$ in the topos $\Psh{\bfC}$ is non-singular w.r.t. ${(\sigma_X :X \rightarrow \Eq \mid X \in \Psh{\bfC})}$ if, and only if, for every figure ${x : \bfC(-, C) \rightarrow X}$, the kernel of $x$ coincides with that of ${f x}$.
\end{example}

We say that a family  ${(\mu_X : X \rightarrow M \mid X \in \calE)}$ satisfies the {\em product property} if the following square commutes
\[\xymatrix{
X\times Y \ar[rrd]_-{\mu_{X \times Y}} \ar[rr]^-{\mu_X \times \mu_Y} && M \times M \ar[d]^-\wedge \\
 & & M
}\]
for every $X$, $Y$ in $\calE$.

\begin{lemma}\label{LemHora4.9} If the family ${\mu}$  satisfies the product property and every split monomorphism is non-singular then 
$\mu$ is a lax cocone.
\end{lemma}
\begin{proof} This is essentially \cite[Lemma~{4.9}]{Hora2024}. If ${f : X \rightarrow Y}$ is a map then then the two triangles below
\[\xymatrix{
X \ar[rrd]_-{\mu_X} \ar[rr]^-{\twopl{id}{f}} && X \times Y  \ar[d]^-{\mu_{X \times Y}} \ar[rr]^-{\mu_X \times \mu_Y} && \ar[lld]^-{\wedge} M \times M \\
   && M
}\]
commute by hypotheses.
\end{proof}

The following concept is somewhat ad-hoc but it is convenient for our purposes.

\begin{definition}\label{DefSingularityGauge}
In a category $\calE$ with finite limits, a {\em singularity gauge} is a meet semilattice ${(\aEq, \wedge, \top)}$ equipped with a family of maps  ${(\sigma_X : X \rightarrow \aEq \mid X \in \calE)}$ satisfying the following:
\begin{enumerate}
\item The product property holds.
\item Monomorphims are non-singular w.r.t. $\sigma$.
\item For every shell ${(\rmB, \beta)}$ on $\calE$,   there exists a unique upper closed subobject ${\aEq_{\beta} \rightarrow \aEq}$ such that, for every $X$ in $\calE$,   ${\sigma_X \beta_X : \rmB X \rightarrow \aEq}$ factors through ${\aEq_{\beta} \rightarrow \aEq}$ and the resulting square
\[\xymatrix{
\rmB X\ar[d]_{\beta_X} \ar@{.>}[r]& \aEq_{\beta} \ar[d]  \\
X \ar[r]_-{\sigma_X} & \aEq
}\]
is a pullback.
\end{enumerate}
\end{definition}

We will give two proofs of the following result.

\begin{theorem}\label{ThmPresheafToposesHaveSingularityGauges} Every presheaf topos has a singularity gauge.
\end{theorem}
\begin{proof}
The quotation from \cite{Lawvere1989} shows that ${(\Eq, \wedge,\top)}$ is a meet-semilattice equipped with a lax-cocone ${(\sigma_X : X \rightarrow \aEq)}$. With the explicit description of $\sigma$ it is easy to check that the product property holds and that monomophisms are non-singular.
As we have already observed in Section~\ref{SecQuotation},  saturated probes on a small category $\bfC$ are in bijective correspondence with upper-closed subobjects of $\Eq$ in $\Psh{\bfC}$. Again, the explicit description of $\sigma_X$ implies that, for a saturated probe $\prP$ in $\bfC$, the pullback of the corresponding subobject of $\Eq$ along $\sigma_X$ coincides with ${\sk_P : \Sk_\prP X \rightarrow X}$. 
\end{proof}

The connection with local state classifiers may be stated as follows.

\begin{theorem}[Hora]\label{ThmHoraLSCimpliesSG} Every local state classifier in a  topos  is a singularity gauge.
\end{theorem}
\begin{proof}
Let $\calE$ be a Grothendieck topos and  let ${(\xi_X : X \rightarrow \Xi  \mid X\in \calE)}$ be a local state classifier therein.
The object $\Xi$ is naturally equipped with a meet semilattice structure by \cite[Proposition~{3.27}]{Hora2024}, and the family ${(\xi_X \mid X \in \calE)}$ is a lax cocone with the product property by \cite[Lemma~{4.9}]{Hora2024}. Monomorphisms are non-singular simply because ${(\xi_X : X \rightarrow \Xi \mid X )}$ is a cocone for ${\calE_m \rightarrow \calE}$.
 Shells are not explicitly considered in \cite{Hora2024} but the argument at the end of \cite[p.296]{Hora2024} makes the role of upper-closedness clear.
\end{proof}

Of course, Theorem~\ref{ThmPresheafToposesHaveSingularityGauges}  follows as a corollary of Theorems~\ref{ThmHoraExistenceOfLSC} and \ref{ThmHoraLSCimpliesSG}.

Naturally, we are interested in a notion of singular map that is absolute in the sense that it is determined by the underlying topos.

\begin{definition}\label{DefNonSingularMap} 
In a topos $\calE$ with local state classifier ${(\xi_X : X \rightarrow \Xi \mid X \in \calE)}$, a morphism ${f : X \rightarrow Y}$ is {\em non-singular} if it is so w.r.t. $\xi$.
\end{definition}

We will discuss a simple example in Section~\ref{SecExample}.
For the moment let us observe the following.

\begin{proposition}\label{PropNonSingularImpliesCartesian} In every topos $\calE$ with a local state classifier, a non-singular map is cartesian w.r.t. any shell on $\calE$.
\end{proposition}
\begin{proof}
Let ${(\sigma_X : X \rightarrow \aEq \mid X \in \calE)}$ be the local state classifier in $\calE$.
Let ${f : X \rightarrow Y}$ be a non-singular map and let ${(\rmB, \beta)}$ be a shell on $\calE$.
By Theorem~\ref{ThmHoraLSCimpliesSG} there exists an upper-closed subobject ${\beta' : \aEq_{\beta} \rightarrow \aEq}$ such that every ${\beta_X}$ is a pullback of $\beta'$ along $\sigma_X$. Then, the outer diagram  and the  right inner square below are pullbacks
\[\xymatrix{
\rmB X \ar[d]_-{\beta_X} \ar[r]^-{\rmB f} & \rmB Y \ar[d]^-{\beta_Y} \ar[r] & \aEq_\beta \ar[d]^-{\beta'} \\
X \ar@(d,d)[rr]_-{\sigma_X} \ar[r]_-f & Y \ar[r]_-{\sigma_Y} & \aEq
}\] 
so, since the left-inner square commutes, it is a pullback by the Pasting Lemma.
\end{proof}

In particular, non-singular maps have discrete fibers when that makes sense.

\section{The petit toposes of non-singular maps}
\label{SecPetitToposes}

Let $\calE$ be a category with finite limits, and let ${(M, \top, \wedge)}$ be a meet semilattice in  $\calE$  equipped with a  lax cocone  ${(\mu_X : X \rightarrow M \mid X \in \calE)}$.

\begin{lemma}\label{LemBasicProp} 
If ${f g}$ is non-singular w.r.t. $\mu$  then so is $g$.
\end{lemma}
\begin{proof}
Let ${g : W \rightarrow X}$ and ${f : X \rightarrow Y}$.
To prove  that $g$ is non-singular calculate
\[ \mu_X g  \leq \mu_Y f g = \mu_W \]
using ${\mu_X \leq \mu_Y f}$ and that $fg$ is non-singular.
\end{proof}

For any object $Y$, let ${\calN_\mu(Y) \rightarrow \calE/Y}$ be the full subcategory of non-singular maps (w.r.t. $\mu$)  with codomain $Y$.
Alternatively, if we let ${\calE_\mu \rightarrow \calE}$ be the subcategory determined by non-singular maps (w.r.t. $\mu$) then
${\calN_\mu(Y) \rightarrow \calE/Y}$  coincides with ${\calE_\mu/Y \rightarrow \calE/Y}$ by Lemma~\ref{LemBasicProp}.

\begin{lemma}\label{LemSingularMapsAreCoreflective} If regular monomorphisms are non-singular then the subcategory ${\calN_\mu(Y) \rightarrow \calE/Y}$ has a right adjoint, and the corresponding counit is a regular mono.
\end{lemma}
\begin{proof}
For a map ${f : X \rightarrow Y}$ we let the fork below
\[\xymatrix{
X' \ar[r]^-{\beta_X} & X \ar@(rd,ld)[rr]_-{\mu_X} \ar[r]^-f  & Y  \ar[r]^-{\mu_Y} & M 
}\]
be an equalizer. The map ${f \beta_X}$ is non-singular because
\[  \mu_Y f \beta_X = \mu_X \beta_X = \mu_{X'} \]
using that the monic $\beta_X$ is  regular  and so, non-singular by hypothesis.

To prove the universal property let ${g : W \rightarrow X}$ be such that ${f g : W \rightarrow Y}$ is non-singular.
That is, ${\mu_Y f g = \mu_W}$. Together with Lemma~\ref{LemBasicProp} we can calculate:
\[ \mu_X g = \mu_W = \mu_Y f g \]
so $g$ factors through ${X' \rightarrow X}$.
\end{proof}

Since identities are non-singular, the subcategory ${\calN_{\mu}(Y) \rightarrow \calE/Y}$ contains the terminal object.
The next result implies that the subcategory is closed under pullbacks.

\begin{lemma}\label{LemPullbacksOfnonSingularMaps} 
If regular monomorphisms are non-singular and the product property holds for $\mu$ then non-singular maps are closed under pullback.
\end{lemma}
\begin{proof}
Assume that ${f : X \rightarrow Y}$ is non-singular and that the square on the left below is a pullback
\[\xymatrix{
P \ar[d]_-h \ar[r]^-k & W \ar[d]^-g &&  P \ar[d]_-{\twopl{h}{k}} \ar[r] & Y \ar[d]^-{\Delta}  \\
X \ar[r]_-f                  & Y                  && X \times W \ar[r]_-{f\times g} & Y \times Y
}\]
so that the square on the right above is also  a  pullback.
Then, since ${\mu_W \leq \mu_Y g}$, we also have ${\mu_W k \leq \mu_Y g k}$ and so
\[ \mu_W k \leq \mu_Y g k = \mu_Y f h = \mu_X h \]
using that $f$ is non-singular. Also, the regular monic  ${\twopl{h}{k}}$ is non-singular by hypothesis, so
\[ \mu_P = \mu_{X \times W} \twopl{h}{k} = \wedge (\mu_X \times \mu_W) \twopl{h}{k} = \wedge \twopl{\mu_X h}{\mu_W k} = \mu_W k \]
using the product property and the previous inequality.
\end{proof}

Combining the above lemmas we obtain the following.

\begin{theorem}\label{ThmPetit} 
Let  $\calE$ be a topos and let  ${(M,  \wedge, \top)}$ be a  meet semilattice therein.
If ${(\mu_X : X \rightarrow M \mid X \in \calE)}$  is a family of maps with the product property and  monos are non-singular (w.r.t. $\mu$) then,  for every object $Y$ in $\calE$,  ${\calN_{\mu}(Y)}$ is a topos and  the inclusion ${\calN_{\mu}(Y) \rightarrow \calE/Y}$  is the inverse image functor of  a hyperconnected geometric morphism ${\calE/Y \rightarrow \calN(Y)}$.
\end{theorem}

Of course, we are mainly interested in the case that the meet semilattice is that of a local state classifier.
Before discussing an example it seems constructive to comment on the relation with the notion of calibration  \cite{Johnstone2012}.
In particular, because it gives an intuition about why we consider the topos ${\calN_{\mu}(Y)}$ a `petit' topos of $Y$.

\begin{definition}\label{DefCalibration} 
A  {\em calibration} in a topos $\calE$ is a subcategory ${\calD \rightarrow \calE}$  satisfying the following properties:
\begin{enumerate}
\item $\calD$ contains all isomorphisms and for any pair of composable maps $f$, $g$  with ${f\in \calD}$, ${f g \in \calD}$ if, and only if, ${g\in \calD}$.
\item $\calD$ is stable under pullback in the sense that, given a pullback as below  %
\[\xymatrix{
P \ar[d]_-h \ar[r]^-k & W \ar[d]^-g   \\
X \ar[r]_-f                  & Y                  
}\]
with ${f\in \calD}$, ${k \in \calD}$.
\item $\calD$ descends along epimorphisms; i.e., in a pullback square as above, if ${k\in \calD}$ and $g$ is epic then ${f\in\calD}$.
\item For every $X$ in $\calE$, the full inclusion ${\calD/X\rightarrow \calE/X}$ has a right adjoint.
\end{enumerate}
\end{definition}

For instance, \cite[Lemma~{2.1}]{Johnstone2012} shows that for a hyperconnected and local ${p : \calE \rightarrow \calS}$, the maps with discrete fibers form a calibration in $\calE$.

In contrast, if we let ${\calD \rightarrow \calE}$ be the subcategory of non-singular maps with respect to some lax cocone then: Lemma~\ref{LemBasicProp} shows that $\calD$  satisfies a strong form of the first axiom of calibrations,
Lemma~\ref{LemPullbacksOfnonSingularMaps} shows that $\calD$  is closed under pullback, 
and Lemma~\ref{LemSingularMapsAreCoreflective}  shows that the relevant subcategories are coreflective.
However, $\calD$  does not descend along epimorphisms in general. We will give a counter-example in \ref{ExFailureOfDescent}.

It is also relevant to observe that the proofs of coreflexivity are quite different.  While the proof Lemma~\ref{LemSingularMapsAreCoreflective} uses an equalizer, the proof of   \cite[Lemma~{2.1}]{Johnstone2012}  cited above uses the fact that toposes are Heyting categories; see also \cite[Proposition~{4.7}]{Menni2022} which makes it clear that the geometric morphism  is not required to be local for the subcategories of maps with discrete fibers to be  coreflective.

So what is the role of localness in the fact that maps with discrete fibers form a calibration?
Localness implies that the direct image functor of the hyperconnected $p$ preserves epimorphisms and, using this, the descent property of calibrations follow. As a side comment we show that the resulting restricted class of hyperconnected maps is `small', regardless of the existence of a local state classifier.

\begin{proposition}\label{PropGoodHyperconnectedAreFew} 
If $\calE$ is an elementary topos then the number of hyperconnected geometric morphisms  from $\calE$ whose direct image functor  preserves epimorphisms is bounded by the number of Lawvere-Tierney topologies in $\calE$.
\end{proposition}
\begin{proof}
If we let $\beta$ be the counit of such a geometric morphism $p$, then \cite[Lemma~{5.5}]{Menni2021} implies that the closure operator ${(\beta_X \Rightarrow -)}$ on subobjects of $X$ is actually universal, so it determines a subtopos of $\calE$.
Let $q$ be another such geometric morphism, with counit $\gamma$,  and assume that $p$ and $q$ determine the same subtopos.
That is, ${(\beta_X \Rightarrow -) = (\gamma_X \Rightarrow -)}$ for each object $X$ of $\calE$.
It easily follows that $\beta_X = \gamma_X$ as subobjects of $X$ and, therefore, ${p = q}$.
\end{proof}

It is natural to wonder if we can improve the results by Rosenthal and by Hora, and also Proposition~\ref{PropGoodHyperconnectedAreFew} by showing that, for any elementary topos $\calE$, the number of hyperconnected maps with domain $\calE$ is small.

\section{Non-singular maps among reflexive graphs}
\label{SecExample}

Let $\Delta$ be the usual category of finite non-empty ordinals and monotone maps between them so that $\Psh{\Delta}$ is the topos of simplicial sets.
Let ${\Delta_1 \rightarrow \Delta}$ be the full subcategory of ordinals with at most two elements so that ${\Psh{\Delta_1}}$ is the topos of reflexive graphs used in \cite{Lawvere2005} to illustrate the distinction between categories `of spaces'  and  `generalized spaces'.

There is another historical curiosity here: the sigularity gauge $\Eq$ in ${\Psh{\Delta}}$ is not ``the loop'' as described by Lawvere in the quotation, but the reflexive graph with a unique point and {\em two} non-identity loops. In more detail, if we denote the objects of ${\Delta_1}$ by ${[0]}$ and ${[1]}$ as usual, where $[0]$ is the terminal object, then the set ${\Eq [1]}$ of edges of the reflexive graph $\Eq$ are the equivalence relations on $[1]$ or, equivalently, the quotients of the associated representable by ${[1]}$, which may be pictured as follows

\medskip
\[\xymatrix{
\Delta_1(- , [1])  = & 
 \cdot \ar@{.>}@(lu,ru) \ar[r] & \cdot \ar@{.>}@(lu,ru) & \ar@{~>}[r] & & \cdot \ar@{.>}@(lu,ru)  \ar@(ld,rd)  &  \ar@{~>}[r] & & \cdot \ar@{.>}@(lu,ru)   & =   1
}\] 
\medskip

\noindent where the dotter arrows represent the distinguished loops, and the wavy arrows represent the quotient maps in the topos.

So the object $\Eq$ in $\Psh{\Delta_1}$ may be pictured as follows.
\medskip
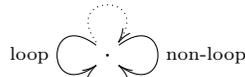
\[\xymatrix{ 
 \ar@(lu,ld)_-{\textnormal{loop}} \cdot \ar@{.>}@(lu,ru)  \ar@(ru,rd)^-{\textnormal{non-loop}}
}\]
and, for each $X$ in $\Psh{\Delta_1}$, the map ${\sigma_X : X \rightarrow \Eq}$ is the unique one that sends non-identity loops in $X$ to `loop', and the rest of non-identity edges to `non-loop'. In other words, very similar to  the case of non-reflexive graphs \cite[Example~{3.9}]{Hora2024}, except for the identity loops.

\begin{remark}\label{RemSimilarity} 
The similarity with \cite[Example~{3.9}]{Hora2024} is, in fact, an important distinction. Notice that the local state classifier in the étendue  of non-reflexive graphs has two points, while that in the  pre-cohesive topos of reflexive graphs has a unique point.  See also Proposition~\ref{PropEqHasAuniquePoint}.
\end{remark}

Using the explicit description of $\sigma$, or Example~\ref{ExNonSingularNaturalTransformations}, we may infer the following.

\begin{corollary}\label{CorNonSingularForDelta1} 
A map in $\Psh{\Delta_1}$ is non-singular if, and only if,  it sends non-identity loops to non-identity loops and non-loops to non-loops.
\end{corollary}

We are now in position to deliver the example promised earlier.

\begin{example}[Non-singular maps in $\Psh{\Delta_1}$ do not form a calibration]\label{ExFailureOfDescent}
If we let 
\[\xymatrix{
A = \Delta_1(- , [1]) \quad = & 
 \cdot \ar@{.>}@(ld,lu) \ar[r] & \cdot \ar@{.>}@(ru,rd) & \ar@{~>}[r]^-q & & \cdot \ar@{.>}@(ld,lu)  \ar@(ru,rd) &   =  L
}\] 
be the obvious quotient in $\Psh{\Delta_1}$ then the kernel pair of  ${q : A \rightarrow L}$ looks as the following pullback
\[\xymatrix{
A + 1 + 1 \ar[d] \ar[r] & A \ar[d]^-q \\
A \ar[r]_-q & L
}\] 
where the top map is non-singular, but the bottom one is not. So the descent property of calibrations does not hold.
\end{example}

Finally we calculate explicitly a couple of petit toposes of non-singular maps and compare them to those calculated in \cite{Lawvere2005}.
Recall that for a reflexive graph $G$, ${\Set(G)}$ denotes the topos of maps with codomain $G$ and  discrete fibers.
It is easy to check that the inclusions ${\calN(1) \rightarrow \Psh{\Delta_1}/1}$ and ${\Set(1) \rightarrow \Psh{\Delta_1}/1}$ coincide with the `discrete' inverse image functor of the canonical geometric morphism ${\Psh{\Delta_1} \rightarrow \Set}$.

For the reflexive graph
\[\xymatrix{
B =  &   \cdot \ar@{.>}@(ld,lu) \ar[r]<+1ex> &  \ar[l]<+1ex> \cdot \ar@{.>}@(ru,rd) 
}\] 
we also have that the inclusions ${\calN(B) \rightarrow \Psh{\Delta_1}/B}$ and ${\Set(B) \rightarrow \Psh{\Delta_1}/B}$ coincide with subcategory of bi-partite graphs. This is because, if ${X \rightarrow B}$ is a map with discrete fibers then $X$ cannot have non-identity loops, so the map must be non-singular.

The next example is motivated by  \cite[Proposition~3]{Lawvere2005} which shows that the étendue ${\Set^{\rightrightarrows}}$ of non-reflexive graphs coincides with the topos ${\Set(L)}$ of maps with discrete fibers over the object $L$ in $\Psh{\Delta_1}$.

\begin{proposition} For the loop $L$ in $\Psh{\Delta_1}$,  ${\calN(L)}$ coincides with the Sierpinski topos. 
\end{proposition}
\begin{proof}
Corollary~\ref{CorNonSingularForDelta1} implies that, for every $X$ in $\Psh{\Delta_1}$, there exists at most one non-singular map  ${X \rightarrow L}$, and $X$ must consist entirely of loops. Also, if the non-horizontal maps below are non-singular
\[\xymatrix{
X \ar[rd] \ar[r]^-f & Y \ar[d] \\
 & L
}\]
then $f$ is non-singular, so it maps non-identity loops to non-identity loops. Altogether, we may forget about the identity loops in ${\calN(L)}$.
\end{proof}

\section{Local state classifiers in pre-cohesive toposes}
\label{SecLSCinPrecohesiveToposes}

The distinction between toposes `of spaces' and `generalized locales' is a powerful guiding principle \cite{Lawvere1989, Lawvere2007, Johnstone2012, Menni2022, Menni2024a} applicable also to the study of local states classifiers.
The main examples of local state classifiers discussed in \cite{Hora2024} are étendues (non-reflexive graphs, localic toposes, Boolean presheaf toposes).
So they are `generalized locales'. 
In previous sections I  emphasized  the case of reflexive graphs which  is an archetypal topos `of spaces'. 

On the side of `generalized locales' we can stress the following two facts.
First, the argument in \cite[Corollary~{5.5}]{Hora2024} implies that if $\calE$ is a topos with local state classifier $\Xi$ then, $\Xi$ is terminal if, and only if, every geometric morphism ${\calE \rightarrow \calS}$ is localic.  In particular, if $\calE$ is a Grothendieck topos, $\Xi$ is terminal if, and only if, $\calE$ is localic.
Second, as Lawvere observes in the quotation of Section~\ref{SecQuotation}, for a small category $\bfC$, the  diagonals ${(\Delta_C \in \Eq C \mid  C \in \bfC)}$ define a point ${1 \rightarrow \Eq}$ in $\Psh{\calC}$ if, and only if, all morphisms in $\bfC$ are monic or, equivalently, ${\Psh{\bfC}}$ is an étendue \cite[Lemma~{C5.2.4}]{elephant}. Hence, on the side of `generalized locales', we are led to wonder about the requirement that the semilattice structure of the local state classider has a bottom element. We will not pursue the idea here, except for the observation below  that it is incompatible with one of the basic properties of toposes `of spaces' 

Recall that a geometric morphism ${p : \calE \rightarrow \calS}$ is {\em local} if the direct image functor has a right adjoint, usually denoted by ${p^! : \calS \rightarrow \calE}$. This is one of the properties of toposes `of spaces' isolated in \cite{Lawvere2005} and is part of the definition of categories of cohesion \cite{Lawvere2007}.

\begin{proposition}\label{PropEqHasAuniquePoint} 
For a small category $\bfC$, if the canonical geometric morphism ${\Psh{\bfC} \rightarrow \Set}$ is local, then the local state classifier  in ${\Psh{\calC}}$ has a unique point.
\end{proposition}
\begin{proof} 
Since we are dealing with presheaf toposes we may, without loss of generality, assume that idempotents split in $\calC$.
The localness hypothesis implies that $\calC$ has a terminal object 1 by  \cite[C3.6.3(b)]{elephant}. 
Then the set of points of $\Eq$ is isomorphic to ${\Eq 1}$, but there is only one equivalence relation on $1$.
\end{proof}

R.~Hora noticed that the same holds for any local ${p  : \calE \rightarrow \Set}$ because ${p_*}$ preserves all colimits.
So a point in the local state classifier must factor through a leg of the universal cocone, and therefore it must coincide with the leg of the cocone determined by the terminal object.

A geometric morphism ${ p : \calE \rightarrow \calS}$ is {\em pre-cohesive} if it is local, hyperconnected, and the inverse image functor $p^*$ has a left adjoint that preserves finite products. The additional left adjoint is denoted by  ${p_! : \calE \rightarrow \calS}$. Intuitively, ${p_! X}$ is the set of connected components of $X$.

Fix a  pre-cohesive  geometric morphism ${ p : \calE \rightarrow \calS}$. 
An object $X$ in $\calE$ is {\em Leibniz} if the canonical map ${p_* X \rightarrow p_! X}$ is an isomorphism. Intuitively,  every connected component has a unique point. For instance, in the topos ${\Psh{\Delta_1}}$, an object is Leibniz if, and only if, it consists of loops only.

It follows from \cite[Theorem~{4.3}]{MarmolejoMenni2021} that  the full subcategory ${\calL \rightarrow \calE}$ of Leibniz objects  is the inverse image functor of a hyperconnected geometric morphism ${s : \calE \rightarrow \calL}$ whose counit will be denoted by ${\lambda_X : s^* (s_* X) \rightarrow X}$ and may be called the {\em Leibniz core} of $X$.

A monomorphism ${u : U \rightarrow X}$ is {\em dense (w.r.t. the centre of $p$)} if it contains the counit ${p^* (p_* X) \rightarrow X}$ of ${p^* \dashv p_*}$. Intuitively, $u$ is dense if it contains all the points of $X$.
Let us say that ${u : U \rightarrow X}$ is {\em lightly dense} if it is dense and it is contained in the Leibniz core of $X$.
Intuitively, it contains all points, but not much more, only infinitesimal information around points.

\begin{definition}\label{DefNotTooDense}
A {\em classifier of lightly dense monomorphisms} is a lightly dense monomorphism ${\Xi_\top \rightarrow \Xi}$ satisfying that, for every lightly dense monomorphism ${u : U \rightarrow X}$, there exists a unique map ${\chi_u : X \rightarrow \Xi}$ with discrete fibers such that ${f \chi_u }$ factors through ${\Xi_\top \rightarrow \Xi}$ and, moreover,  the resulting square
\[\xymatrix{
U \ar[d]_-u \ar[r] & \Xi_{\top} \ar[d] \\
X \ar[r]_-{\chi_u} & \Xi 
}\]
is a pullback.
\end{definition}

Readers may easily identify the category of lightly dense monomorphisms, and pullback squares between them, where the classifier is a terminal object, but to avoid a potential confusion, we stress that lightly dense subobjects are not stable under pullback along maps with discrete fibers. 

We do not know if every pre-cohesive topos has a classifier of lightly dense monomorphisms; but our archetypal  example does.

\begin{proposition}\label{PropNiceCharOfLSCs} 
The pre-cohesive topos ${p : \Psh{\Delta_1} \rightarrow \Set}$ of reflexive graphs has a classifier ${\Xi_\top \rightarrow \Xi}$ of lightly dense monomorphims. In fact, $\Xi$ is the vertex of the  local state classifier, and the component ${\sigma_X :  X \rightarrow \Xi}$ of the universal cocone  is the unique map such that 
\[\xymatrix{
s^*(s_* X) \ar[d]_-{\lambda_X} \ar[r] & \Xi_\top \ar[d] \\
X \ar[r]_-{\sigma_X} & \Xi
} \]
is a pullback. Moreover, ${\Xi_\top \rightarrow \Xi}$ is the image of ${\sigma_{\Xi} : \Xi \rightarrow \Xi}$.
\end{proposition}
\begin{proof}
Recall that we identified  the local state classifier $\Eq$ in $\Psh{\Delta_1}$ as the graph

\medskip
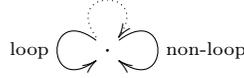
\[\xymatrix{ 
 \ar@(lu,ld)_-{\textnormal{loop}} \cdot \ar@{.>}@(lu,ru)  \ar@(ru,rd)^-{\textnormal{non-loop}}
}\]
and, for each $X$ in $\Psh{\Delta_1}$, the map ${\sigma_X : X \rightarrow \Eq}$ as the unique one that sends non-identity loops in $X$ to `loop', and the rest of non-identity edges to `non-loop'. 
Then, the  image of ${\sigma_{\Eq}}$ is the subobject  ${\ell : \Eq_0 \rightarrow \Eq}$ determined by the edge `loop'.

A monomorphism ${u : U \rightarrow X}$ in $\Psh{\Delta_1}$ is lightly dense if, and only if, it contains all nodes and some of the loops, but no edges between different nodes. So there is a unique map ${\chi_u : X \rightarrow \Eq}$ with discrete fibers such that ${\chi_u^* \ell = u}$.
Indeed, it must send all the non-identity loops in $u$ to `loop' and all the other edges to `non-loop'.
\end{proof}

The phenomenon of Proposition~\ref{PropNiceCharOfLSCs} is reminiscent of the following  basic, but not often stated,  result expressing the terminal object as the colimit of a potentially large diagram.

\begin{proposition}[Freyd-Johnstone?]\label{PropWeirdCharOfFinalObjs}
For any category $\bfC$, terminal objects in $\bfC$ are in bijective correspondence with  colimits of the identity ${\bfC \rightarrow \bfC}$.
\end{proposition}
\begin{proof}
I first saw the result  in  lectures notes of  P. T. Johnstone's  Cambridge Category Theory course.
He claims no originality for it, arguing that Peter Freyd surely knew it when he first proved the Adjoint Functor Theorems.
The proof is not very difficult if you know the statement. It could be a nice exercise, or it may found today  in \cite[Lemma~{3.7.1}]{Riehl2016}. 
See also \cite[Example~{3.12}]{Hora2024}.
%
%
%
\end{proof}

Maybe it is possible strengthen the idea of Proposition~\ref{PropNiceCharOfLSCs}  to characterize local state classifiers in pre-cohesive toposes.


\bibliography{biblio}

\begin{thebibliography}{10}

\bibitem{Hora2024}
R.~Hora.
\newblock Internal parameterization of hyperconnected quotients.
\newblock {\em Theory Appl. Categ.}, 42:263--313, 2024.
\newblock URL: \url{www.tac.mta.ca/tac/volumes/42/11/42-11abs.html#vol42}.

\bibitem{elephant}
P.~T. Johnstone.
\newblock {\em Sketches of an elephant: a topos theory compendium}, volume
  43-44 of {\em Oxford Logic Guides}.
\newblock The Clarendon Press Oxford University Press, New York, 2002.

\bibitem{Johnstone2012}
P.~T. Johnstone.
\newblock Calibrated toposes.
\newblock {\em Bull. Belg. Math. Soc. - Simon Stevin}, 19(5):889--907, 2012.

\bibitem{HoraKamio2024}
Y.~Kamio and R.~Hora.
\newblock A solution to the first {Lawvere}'s problem. {A} {Grothendieck} topos
  that has a proper class many quotient topoi.
\newblock Preprint, {arXiv}:2407.17105 [math.{CT}] (2024), 2024.
\newblock URL: \url{https://arxiv.org/abs/2407.17105}.

\bibitem{LawvereSUNY-Lectures-1975-77}
F.~W. Lawvere.
\newblock Lectures (at {SUNY}, {B}uffalo) 1975-77.
\newblock {U}npublished.

\bibitem{LawvereOpenProblems}
F.~W. Lawvere.
\newblock Open problems in topos theory.
\newblock Talk at the 88th PSSL, Cambridge, England, 4 April 2009.
\newblock URL:
  \url{https://github.com/mattearnshaw/lawvere/blob/master/pdfs/2009-open-problems-in-topos-theory.pdf}.

\bibitem{Lawvere1989}
F.~W. Lawvere.
\newblock Qualitative distinctions between some toposes of generalized graphs.
\newblock Categories in computer science and logic, {Proc}. {AMS}-{IMS}-{SIAM}
  {Jt}. {Summer} {Res}. {Conf}., {Boulder}/{Colo}. 1987, {Contemp}. {Math}. 92,
  261-299, 1989.

\bibitem{Lawvere2005}
F.~W. Lawvere.
\newblock Categories of spaces may not be generalized spaces as exemplified by
  directed graphs.
\newblock {\em Repr. Theory Appl. Categ.}, 2005(9):1--7, 2005.

\bibitem{Lawvere2007}
F.~W. Lawvere.
\newblock Axiomatic cohesion.
\newblock {\em Theory Appl. Categ.}, 19:41--49, 2007.
\newblock URL: \url{https://eudml.org/doc/128088}.

\bibitem{MarmolejoMenni2021}
F.~Marmolejo and M.~Menni.
\newblock The canonical intensive quality of a cohesive topos.
\newblock {\em Theory Appl. Categ.}, 36:250--279, 2021.
\newblock URL: \url{www.tac.mta.ca/tac/volumes/36/9/36-09abs.html}.

\bibitem{Menni2021}
M.~Menni.
\newblock The hyperconnected maps that are local.
\newblock {\em J. Pure Appl. Algebra}, 225(5):15, 2021.
\newblock Id/No 106596.
\newblock \href {https://doi.org/10.1016/j.jpaa.2020.106596}
  {\path{doi:10.1016/j.jpaa.2020.106596}}.

\bibitem{Menni2022}
M.~Menni.
\newblock Maps with discrete fibers and the origin of basepoints.
\newblock {\em Appl. Categ. Struct.}, 30(5):991--1015, 2022.
\newblock \href {https://doi.org/10.1007/s10485-022-09680-2}
  {\path{doi:10.1007/s10485-022-09680-2}}.

\bibitem{Menni2024a}
M.~Menni.
\newblock The {\'e}tendue of a combinatorial space and its dimension.
\newblock {\em Adv. Math.}, 459:21, 2024.
\newblock Id/No 110029.
\newblock \href {https://doi.org/10.1016/j.aim.2024.110029}
  {\path{doi:10.1016/j.aim.2024.110029}}.

\bibitem{Menni2024}
M.~Menni.
\newblock The successive dimension, without elegance.
\newblock {\em Proc. Am. Math. Soc.}, 152(3):1337--1354, 2024.
\newblock \href {https://doi.org/10.1090/proc/16638}
  {\path{doi:10.1090/proc/16638}}.

\bibitem{Riehl2016}
E.~Riehl.
\newblock {\em Category theory in context}.
\newblock Mineola, NY: Dover Publications, 2016.

\bibitem{Rosenthal1982}
K.~I. Rosenthal.
\newblock Quotient systems in {Grothendieck} topoi.
\newblock {\em Cah. Topologie G{\'e}om. Diff{\'e}r. Cat{\'e}goriques},
  23:425--438, 1982.

\end{thebibliography}
\bibliographystyle{plainurl}
\end{document}